\pdfoutput=1
\documentclass{scrartcl}
\author{Elias Pipping\thanks{\emph{E-mail address}: elias.pipping@fu-berlin.de}}
\title{Existence of long-time solutions to dynamic problems of
  viscoelasticity with rate-and-state friction}

\usepackage{nicefrac}

\usepackage{amssymb}
\usepackage{amsmath}
\DeclareMathOperator{\asinh}{asinh}
\newcommand*{\mathd}{{\mathrm d}}

\newcommand*{\onormal}{{\boldsymbol n}}
\newcommand*{\pnstress}{{\bar\sigma_n}}
\newcommand*{\nstress}{{\sigma_n}}
\newcommand*{\tstress}{{\boldsymbol \sigma_t}}
\newcommand*{\stress}{{\boldsymbol\sigma}}
\newcommand*{\strain}{{\boldsymbol\varepsilon}}
\newcommand*{\tenelast}{{\boldsymbol{\mathcal B}}}
\newcommand*{\tenvisco}{{\boldsymbol{\mathcal A}}}

\newcommand*{\bodyforce}{{\boldsymbol b}}

\newcommand*{\testv}{{\boldsymbol v}}
\newcommand*{\testw}{{\boldsymbol w}}

\newcommand*{\velSpace}{V}
\newcommand*{\stateSpace}{X}

\newcommand*{\veca}{{\ddot\vecu}}
\newcommand*{\vecu}{{\boldsymbol u}}
\newcommand*{\vecv}{{\dot\vecu}}

\newcommand*{\vecA}{{\ddot\vecU}}
\newcommand*{\vecU}{{\boldsymbol w}}
\newcommand*{\vecV}{{\dot\vecU}}

\newcommand*{\rmin}[1]{r_{#1}}
\newcommand*{\rref}{{r_\ast}}
\newcommand*{\muref}{{\mu_\ast}}
\newcommand*{\mureg}{{\mu_{\mathrm r}}}
\newcommand*{\mutrunc}{{\mu_{\mathrm t}}}
\newcommand*{\critlen}{{L}}

\newcommand*{\opelast}{{\mathfrak B}}
\newcommand*{\opvisco}{{\mathfrak A}}

\newcommand*{\plog}{\log^+}

\usepackage{xparse}
\usepackage{mathtools}
\DeclarePairedDelimiter\abs{\lvert}{\rvert}
\DeclarePairedDelimiter\paren{\lparen}{\rparen}
\DeclarePairedDelimiter\pnorm{\lvert}{\rvert}
\DeclarePairedDelimiter\norm{\lVert}{\rVert}
\DeclarePairedDelimiter\dualbracket{\langle}{\rangle}
\DeclarePairedDelimiter\curly{\{}{\}}

\usepackage{amsthm}
\usepackage{thmtools}
\declaretheorem{proposition}

\declaretheorem[numberlike=proposition]{theorem}
\declaretheorem[numberlike=proposition,style=definition]{problem}
\declaretheorem[numberlike=proposition]{remark}

\usepackage[T1]{fontenc}

\usepackage{caption}
\usepackage{filecontents}
\usepackage[%
    backend=biber,%
    bibencoding=utf8,%
    doi=true,%
    isbn=false,%
    giveninits=true,%
]{biblatex}
\DeclareFieldFormat{eprint:urn}{%
  URN\addcolon\space
  \ifhyperref
    {\href{http://www.nbn-resolving.org/#1}{\nolinkurl{#1}}}
    {\nolinkurl{#1}}}

\usepackage[autostyle]{csquotes} 
\usepackage[english]{babel}

\newlength\mytemplena
\newlength\mytemplenb
\DeclareDocumentCommand\myalignalign{sm}
{
  \settowidth{\mytemplena}{$\displaystyle #2$}%
  \setlength\mytemplenb{\widthof{$\displaystyle=$}/2}%
  \hskip-\mytemplena%
  \hskip\IfBooleanTF#1{-\mytemplenb}{+\mytemplenb}%
}

\usepackage[inline]{enumitem}
\newlist{assumptions}{enumerate}{10}
\setlist[assumptions]{label*=(A\arabic*)}

\usepackage[hyperfootnotes=false,draft=false]{hyperref}
\usepackage[noabbrev]{cleveref}
\crefname{proposition}{proposition}{propositions}
\crefname{lemma}{lemma}{lemmas}
\crefname{corollary}{corollary}{corollaries}
\crefname{theorem}{theorem}{theorems}
\crefname{problem}{problem}{problems}
\crefname{remark}{remark}{remarks}
\crefname{assumptionsi}{assumption}{assumptions}
\crefname{equation}{}{}

\begin{document}
\maketitle
\abstract{We establish existence of long-time solutions to a dynamic
  problem of bilateral contact between a rigid surface and a
  viscoelastic body, subject to rate-and-state friction. The term
  \emph{rate-and-state friction} is used here to refer to a set of
  functions and equations satisfying conditions which rule out the
  slip law but do cover the ageing law, and thus at least one of the
  rate-and-state friction laws commonly used in the geosciences.

  \medskip


  }
\section{Introduction}
We consider here the dynamic motion of a viscoelastic body
$\Omega \subset \mathbb R^d$ in bilateral contact with a rigid
foundation (on the boundary segment $\Gamma_C$), undergoing
infinitesimal deformation and strain, subject to rate-and-state
friction. To that end, we will derive a weak formulation of the
following problem.
\begin{problem}\label{prob:strong-coupled}
  Find a displacement field $\vecu$ on $\Omega$ of the appropriate
  regularity that satisfies
  \begin{align}
    \label{eq:continuous-displacement-elasticity}
    \stress%
    &= \tenvisco \strain(\vecv) + \tenelast \strain(\vecu)
    &&\text{in $\Omega \times I$}\\
    \label{eq:continuous-displacement-conservation}
    \nabla \cdot \stress + \bodyforce &= \rho \veca
    &&\text{in $\Omega \times I$}\\
    \intertext{with the boundary conditions}
    \label{eq:continuous-displacement-dirichlet}
    \vecv%
    &= 0
    &&\text{on $\Gamma_D \times I$}\\
    \label{eq:continuous-displacement-neumann}
    \stress \onormal%
    & = 0
    &&\text{on $\Gamma_N \times I$}\\
    \label{eq:continuous-displacement-bilateral}
    \vecv \cdot \onormal%
    &= 0
    &&\text{on $\Gamma_C \times I$}\\
    \label{eq:continuous-displacement-friction}
    \myalignalign{-\stress}
    &\left. \myalignalign*{-\tstress}
      \begin{aligned}
        -\tstress
        &= \frac{\mu(\pnorm\vecv,\alpha)\abs\pnstress + C}{\pnorm\vecv} \vecv
        &\text{for $\vecv \ne 0$}\\
        \pnorm\tstress
        &\le \mu(0,\alpha) + C
        &\text{for $\vecv = 0$}
      \end{aligned}\ \right\}%
    &&\text{on $\Gamma_C \times I$}
    \intertext{with prescribed $\vecu(0)$ and $\vecv(0)$ as well as a
       scalar state field $\alpha$ on $\Gamma_C$ that satisfies}
    \label{eq:continuous-state-evolution}
    \dot \alpha + A(\alpha)
    &= f(\pnorm\vecv)
    && \text{on $\Gamma_C \times I$}
  \end{align}
  with prescribed $\alpha(0)$.
\end{problem}
Here, we write $\vecu$ for the displacement, $\bodyforce$ for the body
force, $\stress$ for the stress tensor, and $\tstress$ for its
tangential component where the tangential direction is computed from
the outer normal $\onormal$. Linear Kelvin--Voigt viscoelasticity is
prescribed in~\cref{eq:continuous-displacement-elasticity}, formulated
in terms of the strain tensor $\strain$, a viscosity tensor
$\tenvisco$ and an elasticity tensor $\tenelast$. The friction
law~\cref{eq:continuous-displacement-friction} on $\Gamma_C$ is made
up of the friction coefficient $\mu$, the cohesion $C \ge 0$ and a
prescribed, constant quantity denoted by $\pnstress$, meant to
approximately equal the normal stress $\nstress$. Dirichlet and
Neumann boundary conditions are, furthermore, imposed on the boundary
segments $\Gamma_D$ and $\Gamma_N$, respectively. The mass density is
denoted by $\rho$.

\section{Background}
Rate-and-state friction plays an important role in the modelling of
faults \parencite{doi:10.1016/j.jsg.2010.06.009}, which in turn play
an important role in earthquake nucleation. It expresses frictional
resistance in terms of the sliding velocity or \emph{slip rate}
$\pnorm \vecv$ and an abstract \emph{state} variable $\alpha$. Since
the evolution of this state variable is again governed by the sliding
velocity $\pnorm \vecv$, however, the dependence of the friction
coefficient $\mu$ on the $\pnorm \vecv$ and $\alpha$ should rather be
thought of as a means of depending on $\pnorm \vecv$ in two ways: once
directly, in a monotone fashion, and once indirectly, through
$\alpha$, which reacts less immediately to changes in $\pnorm \vecv$,
but generally in an antitone fashion.

Although laws that go by this name have been derived from experiments
\parencite{doi:10.1029/JB084iB05p02161,%
  doi:10.1029/JB088iB12p10359}, they could just as easily have been
proposed as a regularisation of slip rate dependent friction (in which
the coefficient of friction is a function of the sliding rate only but
the dependence is generally antitone) due to the analytical and
numerical difficulties that such ostensibly simpler stateless laws
present
\parencite{doi:10.1090/qam/1914436}.

The existence and uniqueness of solutions to (weak formulations of)
dynamic problems of viscoelasticity and friction has been thoroughly
studied. Rate-and-state friction falls outside the scope of these
studies, however, because of the variable coupling between the rate
and the state: Neither is typically known. The approach taken in this
work is thus to consider the situation where $\alpha$ is known
a-priori, to then compute $\vecv$ under this assumption (such problems
are covered by the current literature) and to then account for the
actual lack of knowledge of $\alpha$ through a fixed-point iteration.

This work thus parallels earlier work from the author's dissertation
in which the time-discrete setting was
considered \parencite{ThesisPipping}.

\section{Examples}

The following two rate-and-state friction laws are commonly used: the
\emph{ageing law} (also known as \emph{slowness law}), which states
\begin{alignat}{2}
  \label{eq:ageing-law}
  \mu         & = \muref + a\log \frac r{\rref} + b \alpha\text,
  &\qquad
  \dot \alpha &= \frac{\rref e^{-\alpha} - r}\critlen\text,\\
  \intertext{and the \emph{slip law}, which states}
  \label{eq:slip-law}
  \mu         & = \muref + a\log \frac r{\rref} + b \alpha\text,
  &\qquad
  \dot \alpha &= -\frac r\critlen \paren[\big]{\log \frac r\rref + \alpha}\text.
\end{alignat}
When presented in this form, both laws use the same expression for
$\mu$, so that their respective state variables $\alpha$ can be
identified; consequently, the names of these laws are typically used
to refer to the associated state evolution equations only.

The ageing law and the slip law as proposed by
\citeauthor{doi:10.1029/JB084iB05p02161} and
\citeauthor{doi:10.1029/JB088iB12p10359} employ the term
$\log(r/\rref)$, which becomes arbitrarily negative for sliding rates
$r$ close to zero; consequently, we have
\begin{equation*}
  \mu(r,\alpha) \to -\infty
  \qquad
  \text{whenever $r \to 0$}
\end{equation*}
for fixed $\alpha$. They are thus unphysical for sufficiently small
$r$, since they predict a negative coefficient of friction. If we
introduce the quantity
\begin{equation*}
  \rmin\alpha = \rref \exp\paren*{-\frac{\muref + b\alpha}a}\text,
\end{equation*}
this issue becomes even clearer, since now $\mu$ can be written as
\begin{equation}\label{eq:mu-with-rmin}
  \mu(r,\alpha) = a\log \frac r{\rmin\alpha}\text,
\end{equation}
so that $\rmin\alpha$ denotes the rate at which the predicted
coefficient of friction undergoes a sign change. In the literature,
this undesirable behaviour of the
\citeauthor{doi:10.1029/JB084iB05p02161}--%
\citeauthor{doi:10.1029/JB088iB12p10359} laws has been addressed by
means of regularisation
\parencite{doi:10.1016/S0022-5096_01_00042-4}. To be precise, the
logarithm on the right-hand side of~\cref{eq:mu-with-rmin} is replaced
by the nonnegative function $z \mapsto \asinh(z/2)$, yielding the
\emph{regularised law}
\begin{equation}\label{eq:mu-regularised}
  \mureg(r,\alpha) = a\asinh\paren*{\frac r{2\rmin\alpha}}\text.
\end{equation}
A different approach is to trust the original law as much as possible,
and only modify it whenever it predicts a negative coefficient of
friction. The requirement of monotonicity then leads to the
\emph{truncated law}
\begin{equation}\label{eq:mu-truncated}
  \mutrunc(r,\alpha) = a\plog \frac r{\rmin\alpha}
  \qquad \text{with} \qquad
  \plog z = \log \max(1,z)
\end{equation}
Both adjustments clearly guarantee nonnegativity of the friction
coefficient.

In what follows, rather than consider such laws directly, we choose to
work in an abstract setting where friction is described through the
friction coefficient
$\mu \colon \mathbb R^+_0 \times \mathbb R \to \mathbb R^+_0$ and two
functions $A \colon \mathbb R \to \mathbb R$,
$f \colon \mathbb R^+_0 \to \mathbb R$ that govern the state evolution
through the equation
\begin{equation*}
  \dot \alpha + A(\alpha) = f(r)\text.
\end{equation*}
It is immediately clear that the slip law does not fall into this
setting, unfortunately. The ageing law and potentially other laws
of interest, however, do.

\section{Abstract rate-and-state friction}

In working with $\mu$, $A$, and $f$, we find it necessary to make the
following assumptions.
\begin{assumptions}
\item\label{ass:mu-1-monotone} The function $\mu$ is nondecreasing
  and continuous in its first argument.
\item\label{ass:mu-2-lipschitz} The function $\mu$ is uniformly
  Lipschitz in its second argument. In other words, we have
  \begin{equation*}
    \abs{\mu(r, \alpha) - \mu(r, \beta)}
    \le L_\mu \abs{\alpha - \beta}
  \end{equation*}
  for any $\alpha$, $\beta$, and $r \ge 0$.
\item\label{ass:mu-3-upper-bound} The function $\mu$ can be bounded as follows:
  \begin{equation*}
    0 \le \mu(r,\alpha) \le C_\mu(1 + r + \abs{\alpha})
  \end{equation*}
  for any $\alpha$ and $r \ge 0$.%
  \footnote{\Cref{ass:mu-2-lipschitz,ass:mu-3-upper-bound} are not
    independent. Indeed, if we assume the former, the latter reduces
    to requiring $\mu(r, 0) \le C_\mu\paren{1 + r}$.}
\item\label{ass:A-monotone} The function $A$ is nondecreasing and
  continuous.
\item\label{ass:f-lipschitz} The function $f$ is Lipschitz, so that we
  have
  \begin{equation*}
    \abs{f(r) - f(v)} \le L_f\abs{r - v}
  \end{equation*}
  for any $r$ and $v$.
\end{assumptions}
As mentioned earlier, the slip law clearly does not fit into this
framework because of the requirement that $\dot\alpha$ can be written
as a sum of two terms, one of which depends solely on $\alpha$ with
the other depending solely on $r$.

The ageing law, in contrast, satisfies all of the assumptions made
above.
\begin{proposition}
  Consider the ageing law~\cref{eq:ageing-law}, either regularised as
  per~\cref{eq:mu-regularised} or truncated as
  per~\cref{eq:mu-truncated}. Then the resulting law satisfies \cref{%
    ass:mu-1-monotone,%
    ass:mu-2-lipschitz,%
    ass:mu-3-upper-bound,%
    ass:A-monotone,%
    ass:f-lipschitz}.
  \begin{proof}
    That $\mureg$ and $\mutrunc$ satisfy \cref{ass:mu-1-monotone} is
    clear. To show that $\mureg$ satisfies \cref{ass:mu-2-lipschitz},
    it suffices to prove
    \begin{equation*}
      \abs{\mureg(r, \alpha) - \mureg(r, \beta)}
      = a \abs*{
        \asinh\paren*{\frac r{2\rmin{\alpha}}}
        - \asinh\paren*{\frac r{2\rmin{\beta}}}
      }
      \le a \abs*{\log\frac{\rmin{\beta}}{\rmin{\alpha}}}
    \end{equation*}
    for any $\alpha$, $\beta$, and $r \ge 0$, since the right-hand
    side equals $b \cdot \abs{\alpha - \beta}$. For $r = 0$, this is
    immediate; for $r > 0$, it becomes clear once we prove the more
    general claim
    \begin{equation*}
      \abs*{\asinh(x) - \asinh(y)} \le \abs*{\log x - \log y}
    \end{equation*}
    for $x$, $y > 0$. Without loss of generality, assume $x \ge y$, so
    that we need to show
    \begin{equation*}
      \asinh(x) - \asinh(y) \le \log x - \log y\text.
    \end{equation*}
    From the logarithmic representation of the $\asinh$ function, we
    obtain that this is equivalent to
    \begin{equation*}
      \log \frac{x + \sqrt{x^2 + 1}}{y + \sqrt{y^2 + 1}} \le \log \frac xy
    \end{equation*}
    and thus
    \begin{equation*}
      y\sqrt{x^2 + 1} \le x\sqrt{y^2 + 1}
    \end{equation*}
    which is obviously true. For $\mutrunc$, we proceed analogously
    and prove
      \begin{equation*}
        \abs{\mutrunc(r, \alpha) - \mutrunc(r, \beta)}
        = a \abs*{ \plog\frac r{\rmin\alpha} - \plog\frac r{\rmin\beta} }
        \le a \abs*{\log\frac{\rmin\beta}{\rmin\alpha}}\text.
      \end{equation*}
      Again, this is trivially true if $r = 0$. For $r > 0$, we have
      \begin{align*}
        \abs*{ \plog \frac r{\rmin\alpha} - \plog \frac r{\rmin\beta} }
        &= \abs*{
          \log\max \curly*{\frac r{\rmin\alpha},1}
          - \log\max \curly*{\frac r{\rmin\beta},1} }\\
        &= \abs*{
          \max \curly*{\log \frac r{\rmin\alpha},0}
          - \max \curly*{\log \frac r{\rmin\beta},0} }\\
        &\le \abs*{ \log\paren*{\frac r{\rmin\alpha}}
          - \log\paren*{\frac r{\rmin\beta}} }
      \end{align*}
      since $\max\curly{\cdot,0}$ is nonexpansive, so that the claim
      follows. To see that $\mutrunc$ and $\mureg$ satisfy
      \cref{ass:mu-3-upper-bound}, observe only
      \begin{equation*}
        \mutrunc(r,\alpha) = a\plog \frac r{\rmin\alpha}
        \le a\paren*{\plog \frac r\rref + \abs*{\log \frac {\rmin\alpha}\rref}}
        \le a\frac r\rref + \muref + b\abs{\alpha}\text.
      \end{equation*}
      and
      \begin{align*}
        \mureg(r,\alpha)
        &= a\asinh \frac r{2\rmin\alpha}
          = a\log\paren*{
          \frac r{2\rmin\alpha}
          + \sqrt{\paren*{\frac r{2\rmin\alpha}}^2 + 1} }
          \le a\log\paren*{ \frac r{\rmin\alpha} + 1 }\\
        &\le a\log\paren*{2 \max\curly*{1,\frac r{\rmin\alpha}} }
          = a\log2 + \mutrunc(r,\alpha)\text.
      \end{align*}
      Finally, each law clearly satisfies
      \cref{ass:A-monotone,ass:f-lipschitz} with
      \begin{equation*}
        A(\alpha) = -\frac {\rref}\critlen e^{-\alpha}\text,
        \quad
        f(r) = r/\critlen\text,
        \quad \text{and} \quad
        L_f = \critlen\text.\qedhere
      \end{equation*}
  \end{proof}
\end{proposition}
\section{Weak formulation}
Here and in what follows, we will make the following typical
assumptions on the domain $\Omega$, the viscoelastic parameters, the
body force, and the normal stress that we prescribe on the frictional
boundary $\Gamma_C$.
\begin{assumptions}[resume]
\item\label{ass:lipschitz-boundary} The domain $\Omega$ is a bounded
  open subset of $\mathbb R^d$ with a Lipschitz boundary. In
  particular, the $d$-dimensional trace map $\gamma$ is well-defined
  from $H^1(\Omega)^d$ to $L^2(\Gamma)^d$.
\item\label{ass:viscosity-uniformly-elliptic} The viscosity tensor is
  symmetric as well as uniformly bounded from above and below through
  $0 < m_\tenvisco \le M_\tenvisco$, so that
  \begin{equation*}
    m_\tenvisco \norm{\testv}_\velSpace^2
    \le \dualbracket{\opvisco \testv, \testv}
    = \int_\Omega \dualbracket{\tenvisco \strain(\testv),\strain(\testv)}
  \end{equation*}
  and
  \begin{equation*}
    \int_\Omega \dualbracket{\tenvisco \strain(\testv),\strain(\testw)}
    = \dualbracket{\opvisco \testv, \testw}
    \le M_\tenvisco \norm{\testv}_\velSpace \norm{\testw}_\velSpace
  \end{equation*}
  hold for any $\testv$, $\testw \in \velSpace$.
\item\label{ass:elasticity-uniformly-elliptic} The elasticity tensor
  is symmetric as well as uniformly bounded from above and below
  through $0 < m_\tenelast \le M_\tenelast$, so that
  \begin{equation*}
    m_\tenelast \norm{\testv}_\velSpace^2
    \le \dualbracket{\opelast \testv, \testv}
    = \int_\Omega \dualbracket{\tenelast \strain(\testv),\strain(\testv)}
  \end{equation*}
  and
  \begin{equation*}
    \int_\Omega \dualbracket{\tenelast \strain(\testv),\strain(\testw)}
    = \dualbracket{\opelast \testv, \testw}
    \le M_\tenelast \norm{\testv}_\velSpace \norm{\testw}_\velSpace
  \end{equation*}
  hold for any $\testv$, $\testw \in \velSpace$.
\item The body force $\bodyforce$ satisfies
  \begin{equation*}
    \norm{\bodyforce}_{L^2(0,T,\velSpace^*)} < \infty\text.
  \end{equation*}
\item The prescribed normal stress $\pnstress$ satisfies
  \begin{equation*}
    \norm{\pnstress}_{L^\infty(\Gamma_C)} < \infty.
  \end{equation*}
\end{assumptions}

We will work with the spaces
\begin{equation*}
  \velSpace = \curly{
    \testv \in H^1(\Omega)^d \colon
    \text{$\testv = 0$ on $\Gamma_D$, $\testv \cdot \onormal = 0$ on $\Gamma_C$} }
  \quad \text{and} \quad
  H = L^2(\Omega)^d
\end{equation*}
which give rise to the Gelfand triple $V \subset H \subset V^*$, as
well as the space
\begin{equation*}
  \stateSpace = L^2(\Gamma_C)\text.
\end{equation*}
In a standard fashion, by testing
\cref{eq:continuous-displacement-conservation} with functions from
$\velSpace$ at fixed points in time, and putting
\cref{eq:continuous-displacement-elasticity} as well as \cref{%
  eq:continuous-displacement-dirichlet,%
  eq:continuous-displacement-neumann,%
  eq:continuous-displacement-bilateral,%
  eq:continuous-displacement-friction} to use, we obtain the following
weak rate problem.
\begin{problem}\label{prob:weak-rate}
  For given $\alpha \in C(0,T,\stateSpace)$, find
  $\vecu \in L^2(0,T,\velSpace)$ with $\vecv \in L^2(0,T,\velSpace)$
  and $\veca \in L^2(0,T,\velSpace^*)$ such that\footnote{The
    $x$-dependence of each integrand is not made explicit here.}
  \begin{multline}\label{eq:weak-rate}
    \int_\Omega \rho \dualbracket{\veca(t), \testv - \vecv(t)}
    + \int_\Omega \dualbracket{
      \tenvisco \strain(\vecv(t)), \strain(\testv - \vecv(t))}
    + \int_\Omega \dualbracket{
      \tenelast \strain(\vecu(t)), \strain(\testv - \vecv(t))}\\
    + \Phi_\alpha(t,\gamma\testv) - \Phi_\alpha(t,\gamma\vecv(t))
    \ge \int_\Omega \dualbracket{\bodyforce(t), \testv - \vecv(t)}
    \quad \forall \testv \in \velSpace
  \end{multline}
  for almost every $t \in [0,T]$ with prescribed $\vecu(0) = \vecu_0$,
  $\vecv(0) = \vecv_0$ and the friction nonlinearities given by
  \begin{equation*}
    \Phi_\alpha(t,\testv)
    = \int_{\Gamma_C} \varphi_\alpha(t,x,\pnorm{\testv(x)})\,\mathd x
    \quad \text{and} \quad
    \varphi_\alpha(t,x,v) = \int_0^v \mu(r,\alpha(t,x)) \abs\pnstress + C\,\mathd r\text.
  \end{equation*}
\end{problem}
For the state field $\alpha$, meanwhile, we stick to a strong
formulation, requiring the following.
\begin{problem}\label{prob:strong-state}
  For given $\vecv \in L^2(0,T,\velSpace)$, find
  $\alpha \in C(0,T,\stateSpace)$ such that
  \begin{equation*}
    \dot \alpha(t) + A(\alpha(t)) = f(\pnorm{\gamma\vecv(t)})
    \ \text{almost everywhere on $\Gamma_C$}
  \end{equation*}
  for almost every $t \in [0,T]$, with prescribed
  $\alpha(0) = \alpha_0$.
\end{problem}
The reformulation of the coupled \cref{prob:strong-coupled} we will
work with from here on is thus the problem of finding a pair
$(\vecv,\alpha) \in L^2(0,T,\velSpace) \times C(0,T,\stateSpace)$ such
that $\vecv$ solves \cref{prob:weak-rate} with state $\alpha$ and
$\alpha$ solves \cref{prob:strong-state} with rate $\vecv$. To analyse
this problem coupling, we first consider each problem separately
\section{Analysis of the rate problem}
\begin{remark}
  In operator notation, we can also write~\cref{eq:weak-rate} as the
  variational inequality
  \begin{equation}\label{eq:weak-rate-operator}
    \begin{split}
      \rho \dualbracket{\veca(t) + \opvisco \vecv(t) + \opelast \vecu(t) - \bodyforce(t),
        \testv - \vecv(t)}
      &+ \Phi_\alpha(t,\gamma\testv)\\
      &\ge \Phi_\alpha(t,\gamma\vecv(t))
      \quad \forall \testv \in \velSpace
    \end{split}
  \end{equation}
  or the subdifferential inclusion
  \begin{equation}\label{eq:weak-rate-operator-inclusion}
    \bodyforce(t)
    \in \rho \veca(t)
    + \opvisco\vecv(t)
    + \opelast\vecu(t)
    + \gamma^*\partial \Phi_\alpha(t,\cdot)(\gamma\vecv(t))
  \end{equation}
  with $\opvisco$, $\opelast \colon \velSpace \to \velSpace^*$
  given by
  \begin{equation*}
    \opvisco\testv = \int_\Omega
    \dualbracket{\tenvisco \strain(\testv), \strain(\cdot)}
    \quad \text{and} \quad
    \opelast\testv = \int_\Omega
    \dualbracket{\tenelast \strain(\testv), \strain(\cdot)}\text.
  \end{equation*}
\end{remark}
A result on second-order hemivariational inequalities now applies in
particular to our variational setting.
\begin{proposition}\label{prop:rate-problem-solvable}
  \Cref{prob:weak-rate} has a unique solution for any
  $\alpha \in C(0,T,\stateSpace)$, $\vecu_0 \in \velSpace$, and
  $\vecv_0 \in H$.
  \begin{proof}
    For existence of a solution
    see~\textcite[Corollary~12]{doi:10.1016/j.na.2004.11.018}. Uniqueness
    follows in particular from \cref{prop:rate-problem-lipschitz}
    which we prove next.

    A few comments are in order on why Theorem~8 and thus Corollary~12
    from the previously cited work can be applied:
    \Cref{ass:viscosity-uniformly-elliptic,ass:elasticity-uniformly-elliptic}
    make $\opvisco$ and $\opelast$ strongly monotone and symmetric
    bounded linear operators. \Cref{ass:mu-1-monotone}, moreover,
    makes $\varphi_\alpha(t,x,\cdot)$ convex for almost every
    $(t,x) \in [0,T] \times \Gamma_C$, so that the Clarke
    subdifferential of $\varphi_\alpha(t,x,\cdot)$ is actually a
    regular subdifferential. \Cref{ass:mu-3-upper-bound}, finally,
    guarantees
    \begin{equation}\label{eq:phi-derivative-bound}
      \pnorm{\partial \varphi_\alpha(t,x,\cdot)(v)}
      \le C_\mu(1 + \abs v + \abs{\alpha(t,x)}) \abs\pnstress + C\text.
    \end{equation}
    While Theorem~8 in the aforementioned work, as stated, requires
    \cref{eq:phi-derivative-bound} to hold without a $t$- or
    $x$-dependent term, a look at the proof reveals that we are free
    to add any term from $L^2(0,T,\stateSpace)$, and thus in
    particular $\abs{\alpha}$.
  \end{proof}
\end{proposition}
\begin{proposition}\label{prop:rate-problem-lipschitz}
  For two solutions $\vecu$ and $\vecU$ of \cref{prob:weak-rate}
  corresponding to $\alpha$ and $\beta$, respectively, with identical
  initial conditions and $t \in [0,T]$, we have
  \begin{equation*}
    \norm{\vecV-\vecv}_{L^2(0,t,\velSpace)}
    \le
    \sqrt t
    \frac{L_\mu\norm{\gamma}}{m_\tenvisco} \norm{\pnstress}_{L^\infty(\Gamma_C)}
    \norm{\beta - \alpha}_{C(0,t,\stateSpace)}\text.
  \end{equation*}
  In particular, the solution operator $R \colon \alpha \mapsto \vecv$
  is single-valued and Lipschitz with the constant
  \begin{equation*}
    L_R
    = \sqrt T
    \frac{L_\mu\norm{\gamma}}{m_\tenvisco} \norm{\pnstress}_{L^\infty(\Gamma_C)}
  \end{equation*}
  from $C(0,T,\stateSpace)$ to $L^2(0,T,\velSpace)$.
  \begin{proof}
    We test~\cref{eq:weak-rate-operator} for $\vecu$ with $\vecV$ and
    for $\vecU$ with $\vecv$ to obtain
    \begin{align*}
      &\dualbracket{
        \rho(\vecA(s)-\veca(s))
        + \opvisco (\vecV(s)-\vecv(s))
        + \opelast (\vecU(s)-\vecu(s)), \vecV(s) - \vecv(s)}\\
      &\qquad\le \Phi_\alpha(s,\gamma\vecV(s)) - \Phi_\alpha(s,\gamma\vecv(s))
        + \Phi_\beta(s,\gamma\vecv(s)) - \Phi_\beta(s,\gamma\vecV(s))\\
      &\qquad= \int_{\Gamma_C}
        \int_{\pnorm{\gamma\vecv(s)}}^{\pnorm{\gamma\vecV(s)}}
        \paren[\big]{ \mu(r,\alpha) - \mu(r,\beta) } \abs\pnstress\,\mathd r\\
      &\qquad\le L_\mu \int_{\Gamma_C}
        \pnorm{\gamma\vecV(s) - \gamma\vecv(s)}
        \abs{\beta(s) - \alpha(s)} \abs\pnstress\\
      &\qquad\le L_\mu \norm\gamma \norm{\pnstress}_{L^\infty(\Gamma_C)}
        \norm{\vecV(s) - \vecv(s)}_\velSpace
        \norm{\beta(s) - \alpha(s)}_\stateSpace
    \end{align*}
    for almost every $s \in [0,T]$, where the second-to-last estimate
    makes use of \cref{ass:mu-2-lipschitz}. Integrating this
    inequality over the time interval $[0,t] \subset [0,T]$ and
    putting \cref{ass:viscosity-uniformly-elliptic,%
      ass:elasticity-uniformly-elliptic} to use yields
    \begin{multline*}
      \frac \rho 2 \norm{\vecV(t) - \vecv(t)}_H^2
      + m_\tenvisco \norm{\vecV-\vecv}_{L^2(0,t,\velSpace)}^2
      + \frac {m_\tenelast} 2 \norm{\vecV(t) - \vecv(t)}_\velSpace^2\\
      \le
      L_\mu \norm\gamma \norm{\pnstress}_{L^\infty(\Gamma_C)}
      \norm{\vecV - \vecv}_{L^2(0,t,\velSpace)}
      \norm{\beta - \alpha}_{L^2(0,t,\stateSpace)}\text.
    \end{multline*}
    The claim now follows from H\"older's inequality.
  \end{proof}
\end{proposition}

\section{Analysis of the state problem}

In \cref{prob:strong-state}, we view $A$ as an operator on the
function space $\stateSpace$ and obtain a problem that has the
structure of an evolution equation associated with a maximal monotone
operator; in doing so, we do not put the superposition operator
structure of $A$ to use: To solve \cref{prob:strong-state} is to solve
a family of ordinary differential equations at once. In what follows,
we apply the first and second line of thinking, in this order.
\begin{proposition}\label{prop:state-problem-solvable}
  \Cref{prob:strong-state} has a unique solution for any
  $\vecv \in L^2(0,T,\velSpace)$ and $\alpha_0 \in \stateSpace$.
  \begin{proof}
    See for example~\textcite[Theorem~1.3]{doi:10.1007/BF02764629}. We
    remark that the requirement
    \begin{equation*}
      \alpha_0 \in \overline{\operatorname{dom}(A)}
    \end{equation*}
    is automatically fulfilled since we have $L^\infty(\Gamma_C)
    \subset \operatorname{dom}(A)$ and $L^\infty(\Gamma_C)$ is dense
    in $L^1(\Gamma_c)$.
  \end{proof}
\end{proposition}
The solution operator corresponding to
\cref{prop:state-problem-solvable} additionally depends
Lipschitz-continuously on the right-hand side.

\begin{proposition}\label{prop:state-problem-lipschitz}
  For two solutions $\alpha$ and $\beta$ of \cref{prob:strong-state}
  corresponding to $\vecv$ and $\vecV$, respectively, with identical
  initial conditions and $t \in [0,T]$, we have
  \begin{equation}\label{eq:state-lipschitz-pointwise}
    \norm{\alpha(\cdot,x) - \beta(\cdot,x)}_{C(0,t)}
    \le L_f \norm{
      \gamma\vecv(\cdot,x) - \gamma\vecV(\cdot,x)}_{L^1(0,t,\mathbb R^d)}
  \end{equation}
  for almost every $x \in \Gamma_C$ and thus
  \begin{equation}\label{eq:state-lipschitz}
    \norm{\alpha - \beta}_{C(0,T,\stateSpace)}
    \le \sqrt T
    L_f \norm{\gamma\vecv - \gamma\vecV}_{L^2(0,T,\stateSpace^d)}\text.
  \end{equation}
  In particular, the solution operator $S \colon \vecv \mapsto \alpha$ is
  Lipschitz with the constant
  \begin{equation*}
    L_S = \sqrt T \norm{\gamma} L_f
  \end{equation*}
  from $L^2(0,T,\velSpace)$ to $C(0,T,\stateSpace)$.
  \begin{proof}
    For almost every $x \in \Gamma_C$ and $s \in [0,T]$, we have
    \begin{align*}
      \dot\alpha(s,x) + A(\alpha(s,x)) &= f(\pnorm{\gamma\vecv(s,x)})\text,\\
      \dot\beta(s,x) + A(\beta(s,x)) &= f(\pnorm{\gamma\vecV(s,x)})
    \end{align*}
    and thus a pair of evolution equations that have the same
    structure as \cref{prob:strong-state} and are additionally
    one-dimensional. For each such pair we can derive
    \begin{equation*}
      \abs{\alpha(t,x) - \beta(t,x)}
      \le \norm{f(\pnorm{\gamma\vecv(\cdot,x)})
        - f(\pnorm{\gamma\vecV(\cdot,x)})}_{L^1(0,t,\mathbb R^n)}
    \end{equation*}
    for example
    from~\textcite[Theorem~1.2(ii)]{doi:10.1007/BF02764629}). Because
    of \cref{ass:f-lipschitz}, this
    implies~\cref{eq:state-lipschitz-pointwise}. To
    obtain~\cref{eq:state-lipschitz}, we apply H\"older's inequality,
    yielding
    \begin{align*}
      \abs{\alpha(t,x) - \beta(t,x)}
      &\le L_f \norm{
        \gamma\vecv(\cdot,x) - \gamma\vecV(\cdot,x)}_{L^1(0,t,\mathbb R^d)}\\
      &\le \sqrt t L_f \norm{
        \gamma\vecv(\cdot,x) - \gamma\vecV(\cdot,x)}_{L^2(0,t,\mathbb R^d)}
    \end{align*}
    for almost every $(t,x) \in [0,T] \times \Gamma_C$, so that by
    integrating over $\Gamma_C$ we find
    \begin{equation*}
      \norm{\alpha(t,\cdot) - \beta(t,\cdot)}_\stateSpace
      \le \sqrt t L_f \norm{\gamma\vecv - \gamma\vecV}_{L^2(0,t,\stateSpace^d)}\text.
    \end{equation*}
    Since $t \in [0,T]$ was arbitrary, this
    proves~\cref{eq:state-lipschitz}.
  \end{proof}
\end{proposition}

\section{Analysis of the coupled problem}

We first establish short-time existence and uniqueness of a solution.

\begin{proposition}\label{prop:simultaneous-weak-solution}
  For sufficiently small $T > 0$,
  \cref{prob:weak-rate,prob:strong-state} have a unique simultaneous
  solution
  $(\vecv,\alpha) \in L^2(0,T,\velSpace) \times C(0,T,\stateSpace)$
  provided that $\vecu_0 \in \velSpace$, $\vecv_0 \in H$, and
  $\alpha_0 \in \stateSpace$.
  \begin{proof}
    By \cref{prop:rate-problem-lipschitz,%
      prop:state-problem-lipschitz}, the operator
    $R \circ S \colon L^2(0,T,\velSpace) \to L^2(0,T,\velSpace)$ is
    Lipschitz with the constant $L_{RS} = L_RL_S$, which satisfies
    $L_{RS} \to 0$ as $T \to 0$. In particular, the time $T$ can be
    chosen such that we have $L_{RS} < 1$. The claim now follows from
    Banach's fixed point theorem.
  \end{proof}
\end{proposition}

We note that $T$ is not constrained in any way by the values of the
initial data $\vecu_0$, $\vecv_0$ or $\alpha_0$. We can thus extend a
solution provided by \cref{prop:simultaneous-weak-solution} to the
interval $[0,2T]$ by applying the aforementioned proposition
repeatedly: once with the actual initial data to obtain a solution on
the time interval $[0,T]$ and once with the \emph{final data}
resulting from the first application, namely $\vecu(T)$, $\vecv(T)$,
and $\alpha(T)$, to obtain a solution on the interval $[T,2T]$.

That this is indeed possible follows from the embeddings\footnote{By
  $H^1(0,T,\velSpace,\velSpace^*)$ we mean here the space of functions
  $\curly{ \testv \in L^2(0,T,\velSpace) \colon \dot\testv \in
    L^2(0,T,\velSpace^*) }$ equipped with the norm
  $\norm \testv_{H^1(0,T,\velSpace)} =
  \norm{\testv}_{L^2(0,T,\velSpace)} +
  \norm{\dot\testv}_{L^2(0,T,\velSpace^*)}$.}
\begin{equation*}
  \vecu \in H^1(0,T,\velSpace) \subset C(0,T,\velSpace)
  \quad \text{and} \quad
  \vecv \in H^1(0,T,\velSpace,\velSpace^*) \subset C(0,T,H)
\end{equation*}
which give us $\vecu(T) \in \velSpace$ and $\vecv(T) \in H$ in
addition to $\alpha(T) \in \stateSpace$. Since the aforementioned
continuation procedure can be repeated an arbitrary number of times,
we can obtain solutions on $[0,nT]$ for arbitrary $n \in \mathbb N$
and thus intervals of arbitrary size.

\begin{theorem}
  For any $T > 0$, \cref{prob:weak-rate,prob:strong-state} have a
  unique simultaneous solution
  $(\vecv,\alpha) \in L^2(0,T,\velSpace) \times C(0,T,\stateSpace)$
  provided that $\vecu_0 \in \velSpace$, $\vecv_0 \in H$, and
  $\alpha_0 \in \stateSpace$.
\end{theorem}

\section*{Acknowledgements}

This work was supported by the German Research Foundation through
project B01 of the CRC 1114.

\printbibliography
\end{document}